\documentclass{article} 
\usepackage{trees}
\hypersetup{pdftitle={Сколько деревьев в графе},
pdfauthor={Антон Петрунин}}

\begin{document}
\title{Сколько деревьев в графе}
\author{Антон Петрунин}
\date{}
\maketitle

\begin{abstract}
Обсуждается рекурсивная формула для числа остовных деревьев графа.
\end{abstract}

\section{Графы и их остовные деревья}

\begin{wrapfigure}{r}{33 mm}
\vskip-4mm
\begin{tikzpicture}[scale=1.4,
  thick,main node/.style={circle,draw,font=\sffamily\bfseries,minimum size=3mm}]

  \node[main node] (1) at (0,15/6) {$a$};
  \node[main node] (2) at (1,15/6){$b$};
  \node[main node] (11) at (1.5,10/6){$c$};
  \node[main node] (12) at (.5,10/6) {$d$};

  \path[every node/.style={font=\sffamily\small}]
   (11) edge [out=20,in=75,looseness=8] node[above] {} (11)
   (1) edge (2)
   (2) edge[bend left] (12)
   (2) edge[bend right] (12)
   (2) edge (11)
   (11) edge (12)
   (12) edge (12);
\end{tikzpicture}
\end{wrapfigure}

Рассмотрим план шести ежедневных рейсов некоторой авиакомпании между некоторыми парами из аэропортов $a,b,c$ и $d$,
показанный не рисунке.
Для формализации такой и многих других ситуаций в математике используется понятие \emph{граф}.

Граф --- это конечный и не пустой набор \emph{вершин} (в нашем примере вершина графа --- это аэропорт) 
и конечный набор \emph{рёбер}, каждое из которых соединяет пару вершин (в примере ребро --- это рейс авиакомпании).
Пара вершин графа может быть соединена несколькими рёбрами (это может означать, что авиакомпания совершает несколько рейсов в день). 
Также, ребро может соединять вершину с самой собой; в этом случае оно называется \emph{петлёй} (про такое ребро можно думать, как про прогулочный рейс авиакомпании).

Иначе говоря, с математической точки зрения, на диаграмме выше мы видим граф с четырьмя вершинами $a,b,c$ и $d$, 
шестью рёбрами, из них
одна --- петля при вершине $c$, и пара рёбер соединяет $b$ с $d$.
Число рёбер, исходящих из данной вершины, называется её \emph{степенью}, при этом петли считаются дважды;
степени вершин $a,b,c$ и $d$ --- соответственно $1,4,4$ и~$3$.

Изображённый граф является \emph{связным}, то есть, из любой его вершины можно пройти в любую другую пройдя по нескольким его рёбрам.

Предположим, нам требуется сократить число рёбер связного графа, сохранив его связность.
Нетрудно видеть, что это можно сделать тогда и только тогда, когда граф содержит \emph{цикл}.

Цикл --- это маршрут, составленный из различных рёбер, обходящий несколько вершин без повторений 
и возвращающийся в исходную вершину. 
Число рёбер в цикле называется \emph{длиной цикла}.
Например, в нашем графе есть два цикла длины три с вершинами $b$, $c$ и $d$, 
один цикл длины два с вершинами $b$ и $d$,
a также петля при вершине $c$ образует цикл длины один.

Действительно, если из любого цикла связного графа выбросить любое ребро, то граф останется связным.
Более того, если после удаления некоторого ребра граф остался связным, то это ребро принадлежало некоторому циклу --- этот цикл образован самим ребром и кратчайшим путём между его концами в оставшемся графе.

Удаление ребра из цикла можно повторять, пока мы не придём к графу без цикла.
Полученный граф называется остовным деревом исходного графа.

Вообще говоря, связный граф без циклов называется \emph{деревом}.
В~таких графах нет петель, и из любой их вершины в любую другую есть единственный путь по рёбрам без повторений вершин.

На следующем рисунке вы видите все пять различных остовных дерева нашего исходного графа.

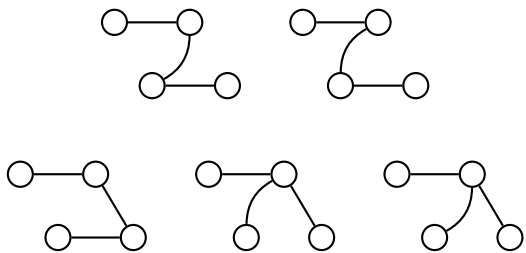
\label{page:5-derev}
\begin{center}
\begin{tikzpicture}[scale=1,
  thick,main node/.style={circle,draw,font=\sffamily\bfseries,minimum size=3mm}]

  \node[main node] (1) at (0,15/6) {};
  \node[main node] (2) at (1,15/6){};
  \node[main node] (11) at (1.5,10/6){};
  \node[main node] (12) at (.5,10/6) {};
  
  \node[main node] (101) at (0+2.5,15/6) {};
  \node[main node] (102) at (1+2.5,15/6){};
  \node[main node] (1011) at (1.5+2.5,10/6){};
  \node[main node] (1012) at (.5+2.5,10/6) {};
  
\node[main node] (201) at (0-1.25,15/6-2) {};
  \node[main node] (202) at (1-1.25,15/6-2){};
  \node[main node] (2011) at (1.5-1.25,10/6-2){};
  \node[main node] (2012) at (.5-1.25,10/6-2) {};
  
  \node[main node] (301) at (0+1.25,15/6-2) {};
  \node[main node] (302) at (1+1.25,15/6-2){};
  \node[main node] (3011) at (1.5+1.25,10/6-2){};
  \node[main node] (3012) at (.5+1.25,10/6-2) {};
  
   \node[main node] (401) at (0+3.75,15/6-2) {};
  \node[main node] (402) at (1+3.75,15/6-2){};
  \node[main node] (4011) at (1.5+3.75,10/6-2){};
  \node[main node] (4012) at (.5+3.75,10/6-2) {};
  
  \path[every node/.style={font=\sffamily\small}]
   (1) edge node{}(2)
   (2) edge[bend left] node{}(12)
   (11) edge node{}(12)

   (101) edge node{}(102)
   (102) edge[bend right] node{}(1012)
   (1011) edge node{}(1012)

   (201) edge node{}(202)
   (202) edge node{}(2011)
   (2011) edge node{}(2012)

   (301) edge node{}(302)
   (302) edge[bend right] node{}(3012)
   (302) edge node{}(3011)
   
   (401) edge node{}(402)
   (402) edge[bend left] node{}(4012)
   (402) edge node{}(4011);
\end{tikzpicture}
\end{center}
Число остовных деревьев в данном графе $\Gamma$ мы будем обозначать $\tau(\Gamma)$.
Например, если $\Gamma$ это граф рассмотренный выше, то $\tau(\Gamma)=5$.

\medskip

Чтобы проверить понимание данных определений, мы советуем решить следующие два стандартных упражнения про деревья.

\begin{thm}{Упражнение}
Докажите, что если дерево имеет хотя бы две вершины, то в нём найдётся вершина степени 1.
\end{thm}

Воспользуйтесь индукцией по числу вершин и предыдущем уп\-ражнением, чтобы доказать следующее.

\begin{thm}{Упражнение}
Докажите, что число рёбер в любом дереве на один меньше числа его вершин.
\end{thm}

В частности, из упражнения следует, что независимо от выбора рёбер,
число удалённых рёбер в процедуре получения остовного дерева из графа описанной выше --- одно и то же.
(Для связного графа $\Gamma$
это число называют его \emph{первым числом Бэтти}; оно обычно  обозначается~$\beta_1(\Gamma)$.)

Ребро связного графа называется \emph{мостом}, если удаление этого ребра из графа делает граф не связным;
такой граф разбивается на два связных графа, называемыми его \emph{островами}.

{
\begin{wrapfigure}{r}{42 mm}
\begin{lpic}[t(-4 mm),b(0 mm),r(0 mm),l(0 mm)]{pics/most(1)}
\lbl[b]{21,8;мост}
\lbl[t]{7,0;остров}
\lbl[t]{34,0;остров}
\end{lpic}
\end{wrapfigure}

\begin{thm}{Упражнение}
Пусть граф $\Gamma$ содержит мост между островами $\Delta_1$ и $\Delta_2$.
Докажите, что
\[\tau(\Gamma)=\tau(\Delta_1)\cdot\tau(\Delta_2).\]
\end{thm}

}

\section{Удаление-плюс-стягивание}\label{sec:deletion+contraction}

Пусть $\rho$ есть ребро в графе $\Gamma$.
Обозначим через $\Gamma\backslash\rho$ граф, полученный из $\Gamma$ удалением ребра $\rho$,
и через $\Gamma/\rho$ граф, полученный из $\Gamma$ стягиванием ребра $\rho$ в точку; смотри рисунок ниже.

Если ребро $\rho$ не является петлёй,  тогда выполняется следующее соотношение
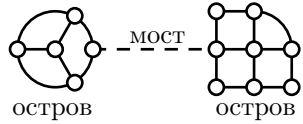
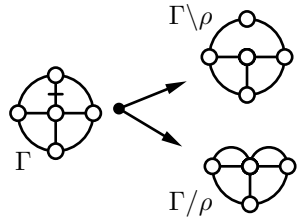
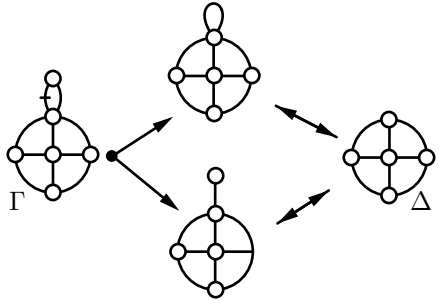
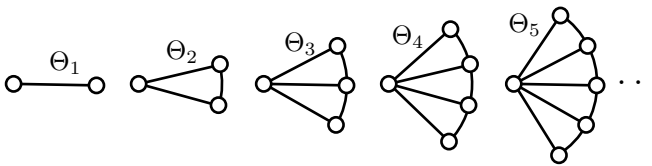
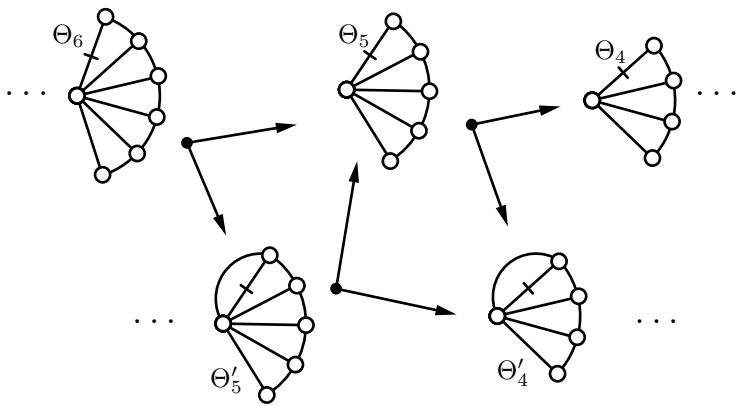
\[\tau(\Gamma)=\tau(\Gamma\backslash\rho)+\tau(\Gamma/\rho),\leqno({*})\]
которое мы будем называть \emph{удаление-плюс-стягивание}.

Действительно, остовные деревья в $\Gamma$ можно разделить на две категории ---
те что содержат ребро $\rho$ и те, что его не содержат.
Для деревьев из первой категории стягивание ребра $\rho$ в точку даёт остовное дерево в $\Gamma/\rho$, а деревья второй категории являются также остовными деревьями в графе  $\Gamma\backslash\rho$.
Более того, оба этих соответствия взаимно однозначны.
Отсюда вытекает формула.

\begin{wrapfigure}{r}{40 mm}
\begin{lpic}[t(-3 mm),b(-3 mm),r(0 mm),l(0 mm)]{pics/osnovnoe-ravenstvo(1)}
\lbl[tr]{4,10;$\Gamma$}
\lbl[br]{28,26;$\Gamma\backslash\rho$}
\lbl[tr]{28,5;$\Gamma/\rho$}
\end{lpic}
\end{wrapfigure}

Например, если $\Gamma$ --- это первый пример и $\rho$ есть ребро между вершинами $b$ и $c$,
тогда первые два остовных дерева на странице \pageref{page:5-derev} соответствуют дереву в $\Gamma\backslash\rho$, а последние два соответствуют дереву в $\Gamma/\rho$.

Формулу $({*})$ удобно записывать схематически, как показано на рисунке.
На графе $\Gamma$ ребро $\rho$, для которого применяется формула, отмечено штрихом. 

Заметим, что никакое остовное дерево не содержит петель.
Поэтому можно удалить все петли из графа, и число его остовных деревьев останется неизменным.
Иначе говоря, для любой петли $\rho$ выполняется равенство 
\[\tau(\Gamma)=\tau(\Gamma\backslash\rho).\]

Из формулы удаление-плюс-стягивание можно вывести несколько других полезных соотношений.
Например, если в графе $\Gamma$ есть концевая вершина $w$ (то есть вершина степени 1), то $w$ и единственное ребро при $w$ можно 
удалить, и в полученном графе $\Gamma\backslash w$ число его остовных графов не изменится, то есть
\[\tau(\Gamma)=\tau(\Gamma\backslash w).\]
Действительно, обозначим через $\rho$ единственное ребро при $w$. 
Заметим, что граф $\Gamma\backslash\rho$ не связен, поскольку вершина $w$ не имеет рёбер, и значит 
$\tau(\Gamma\backslash\rho)=0$.
С другой стороны $\Gamma/\rho=\Gamma\backslash w$, отсюда получаем равенство.

\begin{wrapfigure}{r}{59 mm}
\begin{lpic}[t(-4 mm),b(0 mm),r(0 mm),l(0 mm)]{pics/diagramma(1)}
\lbl[tr]{3,15;$\Gamma$}
\lbl[tl]{54,15;$\Delta$}
\end{lpic}
\end{wrapfigure}

На схемах двусторонняя стрелка ``$\leftrightarrow$'' будет означать, что соответствующие графы имеют то же число остовных деревьев; например, из выведенных тождеств можно вывести следующую диаграмму означающую в частности, что
\[\tau(\Gamma)=2\cdot\tau(\Delta).\]

Равенства, описанные выше, дают алгоритм вычисления $\tau(\Gamma)$.
Действительно, для любого ребра $\rho$, оба графа $\Gamma\backslash\rho$ и $\Gamma/\rho$ имеют меньшее число рёбер.
То есть, удаление-плюс-стягивание сводит нахождение числа остовных деревьев $\Gamma$ к нахождению числа остовных деревьев \emph{более простых} графов.

\section{Деревья в веерах}

Графы следующего вида называются \emph{веерами}; 
веер с $n+1$ вершиной будет обозначаться $\Theta_n$. 

\begin{center}
\begin{lpic}[t(0 mm),b(0 mm),r(0 mm),l(-10 mm)]{pics/veera(1)}
\lbl[br]{12,15;$\Theta_1$}
\lbl[br]{27.5,16.5;$\Theta_2$}
\lbl[br]{43.2,17.5;$\Theta_3$}
\lbl[br]{57.5,18.5;$\Theta_4$}
\lbl[br]{73,20;$\Theta_5$}
\lbl[l]{83,14;{\Large$\dots$}}
\end{lpic}
\end{center}

Применив соотношения, полученные в предыдущей главке, мы можем составить следующую бесконечную схему.
В дополнении к веерам $\Theta_n$ в схеме участвуют их вариации $\Theta_n'$, отличающиеся от $\Theta_n$ дополнительным ребром.
(При возникновении петель или концевых вершин мы их тут же удаляем.)
\begin{center}
\begin{lpic}[t(0 mm),b(0 mm),r(0 mm),l(0 mm)]{pics/veera-skhema(1)}
\lbl[br]{4,49;$\Theta_6$}
\lbl[br]{42,49;$\Theta_5$}
\lbl[br]{76,47;$\Theta_4$}
\lbl[tr]{25,7;$\Theta'_5$}
\lbl[tr]{63,8;$\Theta'_4$}
\lbl[l]{77,13;{\Large$\dots$}}
\lbl[r]{17,13;{\Large$\dots$}}
\lbl[r]{0,43;{\Large$\dots$}}
\lbl[l]{85,43;{\Large$\dots$}}
\end{lpic}
\end{center}

Введём обозначения $a_n=\tau(\Theta_n)$ и $a'_n=\tau(\Theta'_n)$.
Из схемы легко вывести два рекуррентных соотношения:
\begin{align*}
a_{n+1}&=a'_n+a_n,
\\
a'_n&=a_n+a'_{n-1}.
\end{align*}
То есть, в последовательности чисел
\[a_1,a_1',a_2,a_2',a_3\dots\]
каждое следующее является суммой двух предыдущих.

Напомним, что последовательность \emph{чисел Фибоначчи} $F_n$ задаётся тем же соотношением 
$F_{n+1}=F_n+F_{n-1}$ с $F_1=F_2=1$.
Она начинается следующим образом
\[1,1,2,3,5,8,13,\dots\]

Далее заметим, что $\Theta_1$ --- это две вершины соединённые единственным ребром,
а $\Theta'_1$ --- это две вершины соединённые двойным ребром.
Отсюда $a_1=1=F_2$ и $a_1'=2=F_3$ и значит 
\[a_n=F_{2\cdot n}\]
для любого $n$.

Можно также вывести рекуррентное соотношение на $a_n$, без использования $a_n'$:
\begin{align*}
a_{n+1}&=a_n'+a_n=
\\
&=2\cdot a_n+a'_{n-1}=
\\
&=3\cdot a_n-a_{n-1}.
\end{align*}

Последовательности, для которых выполняется подобное соотношение, называются \emph{линейными рекуррентными последовательностями}.
Для таких последовательностей можно найти формулу для её общего члена.
В нашем случае \[a_n=\tfrac1{\sqrt{5}}\cdot
\left(
(\tfrac{3+\sqrt{5}}2)^n-(\tfrac{3-\sqrt{5}}2)^n
\right).\]
Этот процесс подробно описан в книжке \cite{markushevich}, которую мы рекомендуем читателю.

\medskip

Для закрепления материала мы советуем разобрать следующие упражнения:

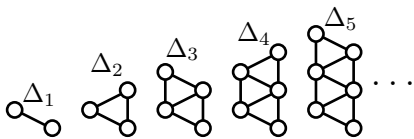
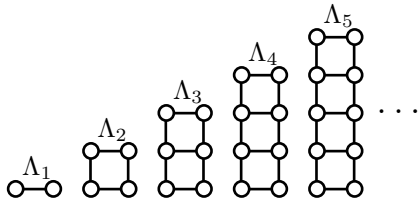
\begin{thm}{Упражнение}
Рассмотрим последовательность графов $\Delta_n$ следующего типа:
\begin{center}
\begin{lpic}[t(1 mm),b(0 mm),r(0 mm),l(0 mm)]{pics/a-extra(1)}
\lbl[b]{5.5,5;$\Delta_1$}
\lbl[b]{14.5,9;$\Delta_2$}
\lbl[b]{24.5,11;$\Delta_3$}
\lbl[b]{34,13;$\Delta_4$}
\lbl[b]{45.5,15;$\Delta_5$}
\lbl[l]{49.5,8;{\Large$\dots$}}
\end{lpic}
\end{center}

Докажите, что $\tau(\Delta_n)=\tau(\Theta_n)$, и значит $\tau(\Delta_n)=F_{2\cdot n}$ для любого $n$. 

\end{thm}

\begin{thm}{Упражнение}
Пусть $b_n$ обозначает число остовных деревьев в \emph{лестнице} $\Lambda_n$ с $n$ ступеньками, то есть в графе следующего типа.

\begin{center}
\begin{lpic}[t(1 mm),b(0 mm),r(0 mm),l(0 mm)]{pics/a1-4(1)}
\lbl[b]{5.5,4;$\Lambda_1$}
\lbl[b]{15.5,9;$\Lambda_2$}
\lbl[b]{25.5,14;$\Lambda_3$}
\lbl[b]{35.5,19;$\Lambda_4$}
\lbl[b]{45.5,24;$\Lambda_5$}
\lbl[l]{50.5,13;{\Large$\dots$}}
\end{lpic}
\end{center}

Воспользуйтесь разработанным методом и докажите, что последовательность $b_n$ удовлетворяет следующему соотношению 
\[b_{n+1}=4\cdot b_n-b_{n-1}.\]

\end{thm}

\begin{wrapfigure}{r}{20 mm}
\begin{lpic}[t(-0 mm),b(0 mm),r(0 mm),l(0 mm)]{pics/lestnitza-shtrih(1)}
\lbl[b]{5,23.5;$\Lambda_3'$}
\lbl[b]{15,21;$\Lambda_3''$}
\end{lpic}
\end{wrapfigure}

Заметим, что $b_1=1$ и $b_2=4$; отсюда можно быстро посчитать первые члены последовательности:
\[1,4,15,56,209,780,2911,\dots \]

Вам потребуются ещё пара последовательностей графов, 
$\Lambda_n'$ и $\Lambda_n''$ показанных на рисунке.

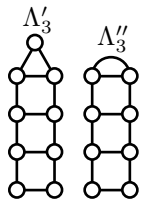
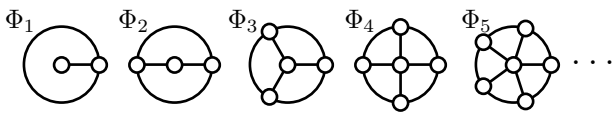
\begin{thm}{Продвинутое упражнение}
Колёсами назыаются графы  $\Phi_n$ следующего типа:
\begin{center}
\begin{lpic}[t(1 mm),b(0 mm),r(0 mm),l(0 mm)]{pics/kolesa(1)}
\lbl[br]{3,11;$\Phi_1$}
\lbl[b]{16,11;$\Phi_2$}
\lbl[b]{30.5,11;$\Phi_3$}
\lbl[b]{46,11;$\Phi_4$}
\lbl[b]{61.5,11;$\Phi_5$}
\lbl[l]{74,7;{\Large$\dots$}}
\end{lpic}
\end{center}
Докажите, что последовательность $c_n=\tau(\Phi_n)$ удовлетворяет следующей рекуррентной формуле:
\[c_{n+1}=4\cdot c_n-4\cdot c_{n-1}+c_{n-2}.\]

\end{thm}

При оказывается проще вывести соотношение 
\[c_{n+1}-a_{n+1}= c_n+a_n.\]
Зная, что $c_1=1$, $c_2=5$, последнее равенство позволяет получить выражение для $c_n$
\[c_n=F_{2\cdot n+1}+F_{2\cdot n-1}-2=L_{2\cdot n}-2;\]
здесь $L_n=F_{n-1}+F_{n+1}$; эта последовательность называются \emph{числами Люка}.
Как и числа Фибоначчи, числа Люка часто проявляются в различных комбинаторных задачах.

\section{Графы как матрицы}

{

\begin{wrapfigure}[4]{r}{28 mm}
\vskip -5mm
\begin{tikzpicture}[scale=1.4,
  thick,main node/.style={circle,draw,font=\sffamily\bfseries,minimum size=3mm}]

  \node[main node] (1) at (0,15/6) {$1$};
  \node[main node] (2) at (1,15/6){$2$};
  \node[main node] (11) at (1.5,10/6){$3$};
  \node[main node] (12) at (.5,10/6) {$4$};

  \path[every node/.style={font=\sffamily\small}]
  
   (1) edge node{}(2)
   (2) edge[bend left] node{}(12)
   (2) edge[bend right] node{}(12)
   (2) edge node{}(11)
   (11) edge node{}(12)
   (12) edge node{}(12);
\end{tikzpicture}
\end{wrapfigure}

Пронумеруем вершины графа $\Gamma$ числами от $1$ до $n$.
Тогда граф $\Gamma$ можно полностью описать таблицей $n\times n$, поставив в клетки на пересечении $i$-ой строки и $j$-ого столбца число рёбер, соединяющих $i$-тую вершину графа с и $j$-ой.

}

Полученная таблица $A=A_\Gamma$ называется \emph{матрицей смежности графа}.
Она симметрична, то есть отражение в главной диагонали переставляет равные числа.
Для графа на рисунке получим 

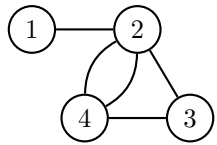
\[A=\left(
\begin{matrix}
0&1&0&0
\\
1&0&1&2
\\
0&1&0&1
\\
0&2&1&0
\end{matrix}
\right)\]

Напомним, что нас интересует число остовных деревьев в графе.
Поскольку наличие петель не влияет на это число, мы можем считать, что граф не содержит петель.
В этом случае на главной диагонали $A$ стоят нули.

Поскольку матрица смежности полностью описывает исходный граф, 
число остовных деревьев графа в принципе можно вычислить по его матрице смежности.

Оказывается, это довольно легко сделать.
Для этого, из построенной $n\times n$-матрицы  $A=A_\Gamma$ для графа $\Gamma$ нужно построить так называемый \emph{минор Кирхгофа} --- $(n-1)\z\times(n-1)$-матрицу определённого типа. 

\begin{enumerate}
\item Обратим знаки всех компонент матрицы $A$ и заменим нули на её диагонали степенями соответствующих вершин. 
В полученной матрице $A'$ сумма чисел в каждой строке и в каждом столбце равна нулю. 
\item Удалим из матрицы $A'$ последнюю строку и последний столбец;
это и есть \emph{минор Кирхгофа} графа $\Gamma$, далее 
обозначаемый $K=K_\Gamma$.
\end{enumerate}

Для примера выше имеем
\[A'=\left(
\begin{matrix}
1&-1&0&0
\\
-1&4&-1&-2
\\
0&-1&2&-1
\\
0&-2&-1&3
\end{matrix}
\right),
\quad 
K=\left(
\begin{matrix}
1&-1&0
\\
-1&4&-1
\\
0&-1&2
\end{matrix}
\right)\]

Заметим, что по минору Кирхгофа $K$ можно восстановить $A'$, а значит и матрицу $A$, и граф $\Gamma$, если у него не было петель.
Действительно, чтобы получить $A'$ нужно добавить строку и столбец к минору Кирхгофа $K$ и заполнить их числами так, чтобы сумма в каждой строке и столбце была нулевой;
это можно сделать единственным образом.

Таким образом, по минору Кирхгофа $K$ можно найти и число деревьев исходного графа $\tau(\Gamma)$.
Мы обозначим это число как $\delta(K)$, подчёркивая этим, что $\delta(K)=\tau(\Gamma)$ есть функция от $K$.
Для $\delta(K)$ верны следующие соотношения:
\begin{enumerate}
\item\label{1-delta} Если в миноре Кирхгофа $K$ переставить пару строк и так же пару столбцов с теми же номерами,
тогда для полученного минора Кирхгофа $K'$ выполняется 
\[\delta(K)=\delta(K').\]
Действительно, $K'$ описывает тот же граф $\Gamma$ с другой нумерацией вершин.
\item\label{2-delta}
Если сумма компонент первой строки в $K$ положительна, то
\[\delta(K)=\delta(K^{\circ})+\delta(K^{\bullet})\]
где $K^{\circ}$ обозначает матрицу $K$, в которой от углового компоненты отняли 1, а $K^{\bullet}$ обозначает $(n-2)\times(n-2)$ матрицу, полученную из $K$ удалением первого столбца и первой строки.
Это равенство следует из формулы удаление-плюс-стягивание $({*})$ поскольку \[K^{\circ}=K_{\Gamma\backslash\rho}\quad\text{и}\quad K^{\bullet}=K_{\Gamma/\rho}.\]
\item Если сумма чисел в каждой строке $K$ равна $0$ то $\delta(K)=0$. 
Действительно, в этом случае ни одна из вершин исходного графа не связана с последней вершиной, а, значит, граф несвязен и не имеет остовных деревьев.

\item\label{4-delta} Если $K$ есть единичная матрица, то есть

\[
K=
\left(
\begin{matrix}
1&0&\cdots&0
\\
0&1&\ddots&\vdots
\\
\vdots&\ddots&\ddots&0
\\
0&\cdots&0&1
\end{matrix}
\right),
\]
тогда $\delta(K)=1$.
Действительно, в этом случае все вершины графа $\Gamma$ с номерами $1,
\dots,n-1$ соединены с вершиной номер $n$ единственным ребром.
Значит граф $\Gamma$ сам является деревом и $\tau(\Gamma)=1$.

\item Свойства \ref{1-delta}---\ref{4-delta} полностью описывают $\delta(K)$ для любой матрицы Кирхгофа.
Последнее следует из существования алгоритма для нахождения $\tau(\Gamma)$ описанный в главе  \ref{sec:deletion+contraction}.
\end{enumerate}

Читатель, знакомый с понятием определитель%
\footnote{Для читателей, незнакомых с понятием определитель, написана следующая главка.},
непременно заметит, что
тем же свойствам \ref{1-delta}---\ref{4-delta} обладает и \emph{определитель матрицы} $K$, далее обозначаемый $|K|$.
Поскольку эти свойства полностью определяют $\delta(K)$, мы получаем равенство
\[\delta(K)=| K|.\]

Отсюда следует так называемая \emph{матричная формула}
\[\tau(\Gamma)=| K_\Gamma|.\leqno({*}{*})\]

Идея доказательства поучительна ---  равенство двух раз\-нo\-шёр\-стно определённых числел
$\tau(\Gamma)$ и $|K_\Gamma|$ следует из общих свойств этих чисел которые полностью определяют $\tau(\Gamma)$. 
Эту идею можно рассматривать как обобщение метода математической индукции;
она имеет множество других приложений, но, к сожалению, более простые содержательные примеры нам не известны.

\section{Ликбез на тему «определители»}

\emph{Квадратная матрица} это таблица $n{\times}n$, заполненная числами, называемыми её \emph{компонентами}.
Определитель $| M|$ матрицы $M$ это многочлен от $n^2$ её компонент,
который удовлетворяет следующим условиям:
\begin{enumerate}
 \item\label{1} Определитель единичной матрицы равен 1; то есть,
\[
\left|
\begin{matrix}
1&0&\cdots&0
\\
0&1&\ddots&\vdots
\\
\vdots&\ddots&\ddots&0
\\
0&\cdots&0&1
\end{matrix}
\right|=1.
\]
\item\label{2} Если каждую компоненту одной из строк матрицы $M$ умножить на число $\lambda$, то определитель полученной матрицы $M'$ умножится на то же число, то есть
\[|M'|=\lambda\cdot |M|.\]
\item\label{3} Если одну из строк матрицы $M$ почленно прибавить к другой строке (или отнять от неё), то определитель полученной матрицы $M'$ не изменится, то есть
\[|M'|= |M|.\]
\end{enumerate}

Эти три свойства однозначно определяют определитель.
Мы примем это утверждение без доказательства; оно неочевидное, но и несложное, 
кроме того, рано или поздно вам придётся его выучить.
То же относится к тождеству 
\[|M|=|M^{\circ}|+|M^{\bullet}|\]
соответствующего свойству \ref{2-delta} для функции $\delta$ в предыдущей главе;
оно следует из так называемого \emph{разложения определителя по строке}.

\begin{thm}{Упражнение}
Выведите следующее свойство из трёх перечисленных выше.
\end{thm}

\begin{enumerate}[resume]
 \item 
Если две строки матрицы $M$ поменять местами, то определитель полученной матрицы $M'$ поменяет знак; то есть,
\[|M'|=-|M|.\]
\end{enumerate}

Для определителя матрицы $n\times n$ можно выписать явную формулу из $n!$ слагаемых. 
Например, 
\[
a_1{\cdot} b_2{\cdot} c_3+a_2{\cdot} b_3{\cdot} c_1+a_3{\cdot} b_1{\cdot} c_2-a_3{\cdot} b_2{\cdot} c_1-a_2{\cdot} b_1{\cdot} c_3-a_1{\cdot} b_3{\cdot} c_2\]
есть определитель матрицы
\[\left(
\begin{matrix}
a_1&b_1&c_1
\\
a_2&b_2&c_2
\\
a_3&b_3&c_3
\end{matrix}
\right).\]
Однако, свойства, описанные выше, дают более удобный и быстрый способ вычисления определителя, особенно при больших значениях $n$.
Мы разберём этот способ на одном примере, который нам пригодится позже:
\[\left|
\begin{matrix}
4&-1&-1&-1
\\
-1&4&-1&-1
\\
-1&-1&4&-1
\\
-1&-1&-1&4
\end{matrix}
\right|
\ 
=
\ 
\left|
\begin{matrix}
1&1&1&1
\\
-1&4&-1&-1
\\
-1&-1&4&-1
\\
-1&-1&-1&4
\end{matrix}
\right|
\ 
=
\ 
\left|
\begin{matrix}
1&1&1&1
\\
0&5&0&0
\\
0&0&5&0
\\
0&0&0&5
\end{matrix}
\right|
\ 
=
\ 
\]
\[
=\ 
5^3\cdot
\left|
\begin{matrix}
1&1&1&1
\\
0&1&0&0
\\
0&0&1&0
\\
0&0&0&1
\end{matrix}
\right|\ =\ 
5^3\cdot\left|
\begin{matrix}
1&0&0&0
\\
0&1&0&0
\\
0&0&1&0
\\
0&0&0&1
\end{matrix}
\right|
\ =\ 5^3.\]
Первое равенство следует из свойства \ref{3} --- мы по очереди прибавили к первой строке все строки со второй до последней. 
Далее мы прибавили первую ко всем остальным, и к полученной матрице применили свойство \ref{2} три раза.
Последние два равенства получаются отниманием от первой строки всех остальных, и последующим применением свойства \ref{1}.

\section{Формула Кэли}

\emph{Полным графом} называется граф в котором любая пара различных вершин соединена единственным ребром;
полный граф с $n$ вершинами будет обозначаться $\Pi_n$.

\begin{center}
\begin{lpic}[t(1 mm),b(0 mm),r(0 mm),l(0 mm)]{pics/polnyj-graf(1)}
\lbl[b]{2.5,8;$\Pi_1$}
\lbl[b]{17,8;$\Pi_2$}
\lbl[bl]{32.3,9.3;$\Pi_3$}
\lbl[bl]{49,10;$\Pi_4$}
\lbl[bl]{65,10;$\Pi_5$}
\lbl[l]{72,6.5;{\Large$\dots$}}
\end{lpic}
\end{center}

Легко видеть, что каждая вершина графа $\Pi_n$ имеет степень $n-1$
и значит минор Кирхгофа $K$ в матричной формуле $({*}{*})$ для $\Pi_n$ есть $(n-1)\times (n-1)$ матрица следующего вида:
\[
K=\left(
\begin{matrix}
n{-}1&-1&\cdots&-1
\\
-1&n{-}1&\ddots&\vdots
\\
\vdots&\ddots&\ddots&-1
\\
-1&\cdots&-1&n{-}1
\end{matrix}
\right).
\]

Прямое обобщение вычисления определителя данного в предыдущей главке даёт, что
\[|K|=n^{n-2}\]
и значит
\[\tau(\Pi_n)=n^{n-2}\]
для любого $n$.

\begin{wrapfigure}[3]{r}{15 mm}
\begin{lpic}[t(-30 mm),b(0 mm),r(0 mm),l(0 mm)]{pics/Pi34(1)}
\end{lpic}
\end{wrapfigure}

Это равенство называется \emph{формулой Кэли}.

Пусть $\Pi_{m,n}$ есть \emph{полный двудольный граф}, то есть граф вершины которого разделены на две доли с $m$ и $n$ вершинами соответственно и у которого любая вершина первой доли соединена со всеми вершинами второй доли.
На рисунке выше изображён полный двудольный граф $\Pi_{4,3}$, вершины каждой доли отмечены своим цветом.

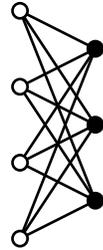
\begin{thm}{Упражнение}
Воспользуйтесь матичной формулой, чтобы доказать равенство
\[\tau(\Pi_{m,n})=m^{n-1}\cdot n^{m-1}.\]

\end{thm}

\section{Заключительные замечания}

Подсчёт числа остовных деревьев связан с расчётом электрических цепей.

Предположим граф $\Gamma$ описывает электрическую цепь;
каждое ребро это единичное сопротивление в омах.
Источник питания подключён к вершинам $a$ и $b$ и $I$ есть общий ток в цепи составляет в амперах.
Нам надо посчитать ток через одно из сопротивлений.

Пусть $\rho$ есть ребро $\Gamma$ с выбранным направлением.
Заметим, что все остовные деревья $\Gamma$ можно разделить на три типа,
(1) те в которых на пути из $a$ в $b$ ребро $\rho$ появляется с положительной ориентацией,
(2) те в которых на пути из $a$ в $b$ ребро $\rho$ появляется с отрицательной ориентацией,
(3) те в которых на пути из $a$ в $b$ ребро $\rho$ не появляется.
Обозначим через $\tau_+$, $\tau_-$ и $\tau_0$ число деревьев в этих трёх категориях.
Очевидно, что $\tau(\Gamma)=\tau_++\tau_-+\tau_0$.

Силу тока $I_\rho$ вдоль $\rho$ можно вычислить по следующей формуле:
\[I_\rho=\frac{\tau_+-\tau_-}{\tau(\Gamma)}\cdot I.\]
Доказательство получается прямой проверкой правил Кирхгофа для значений токов полученных по этой формуле для всех рёбер в~$\Gamma$.

Есть множество других приложений правил Кирхгофа в теории графов. 
Например в \cite{levi}, они используются в физическом доказательстве формулы Эйлера
\[\text{В}-\text{Р}+\text{Г}=2,\]
где $\text{В}$, $\text{Р}$ и $\text{Г}$ обозначает число вершин, рёбер и граней многогранника.

Формула \emph{удаление-плюс-стягивание} применялась при решении так называемой \emph{задачи о квадрировании квадрата}.
История этой задачи и её замечательное решение обсудаются в книжках \cite{yaglom} и \cite[Глава 32]{gardner};
идея этого решения также основана на оригинальном использовании электрических цепей.

Для числа раскрасок вершин графа выполняется аналогичная формула, \emph{удаление-минус-стягивание}.
А именно, пусть $\chi(\Gamma,k)$ обозначает число раскрасок графа $\Gamma$ в $k$ цветов при которых концы каждого ребра покрашены в разные цвета.
Тогда выполняется соотношение
\[\chi(\Gamma,k)=\chi(\Gamma\backslash\rho,k)-\chi(\Gamma/\rho,k).\]

Действительно, допустимые раскраски графа $\Gamma\backslash\rho$ можно разбить на две категории: (1) те в которых концы ребра $\rho$ покрашены в разные цвета --- такие остаются допустимыми в графе $\Gamma$ и (2) те в которых концы ребра $\rho$ покрашены в один цвет --- каждой такой раскраске соответствует единственная раскраска $\Gamma/\rho$.
Отсюда формула.

\begin{thm}{Упреажнение}
Докажите что для любого графа $\Gamma$, функция $\chi(\Gamma,k)$ является многочленом с целыми коэффициентами от $k$;
он называется \emph{хроматическим многочленом} графа $\Gamma$.

(Подсказка: воспользуйтесь индукцией по общему числу рёбер и вершин графа и формулой удаление-минус-стягивание.)
\end{thm}

Приведённый нами вывод рекуррентных формул для чисел ос\-тов\-ных деревьев в веерах лестницах и колёсах даётся в \cite{haghighi-bibak};
эта задача обсуждается также в классической книжке \cite{knut}.
Наиболее известным доказательством формулы Кэли является так называемый \emph{код Прюфера} --- способ однозначного кодирования дерева упорядоченной последовательностью из его вершин; нескольким другим доказательствам посвящена Глава 30  в \cite{aigner-ziegler}.

\end{document}